\documentclass[12pt]{amsart}
\usepackage{amsmath}
\usepackage{amssymb}
\usepackage{latexsym}
\usepackage{amscd}
\usepackage{citesort}

\newdimen\AAdi%
\newbox\AAbo%
\def\AAk#1#2{\s_etbox\AAbo=\hbox{#2}\AAdi=\wd\AAbo\kern#1\AAdi{}}%
\def\AAr#1#2#3{\s_etbox\AAbo=\hbox{#2}\AAdi=\ht\AAbo\raise#1\AAdi\hbox{#3}}%
\font\tenmsb=msbm10 at 12pt \font\sevenmsb=msbm7 at 8pt
\font\fivemsb=msbm5 at 6pt
\newfam\msbfam
\textfont\msbfam=\tenmsb \scriptfont\msbfam=\sevenmsb
\scriptscriptfont\msbfam=\fivemsb

\newtheorem{thm}{Theorem}[section]
\newtheorem{lem}{Lemma}[section]

\newtheorem{rem}{Remark}[section]
\newtheorem{pro}{Proposition}[section]

\newcommand{\ba}{\begin{array}}
\newcommand{\ea}{\end{array}}

\newcommand{\Section}[2]{\setcounter{equation}{0}
\allowdisplaybreaks
\section[#1]{#2}}

\def\n{\nabla}

\def\ir#1{\mathbb R^{#1}}

\def\f#1#2{\frac{#1}{#2}}

\def\grs#1#2{\bold G_{#1,#2}}

\def\mc#1{\mathcal{#1}}

\def\pd#1#2{\frac {\partial #1}{\partial #2}}

\def\td{\tilde}

\def\a{\alpha}
\def\be{\beta}

\def\p#1{\partial #1}

\def\de{\delta}
\def\De{\Delta}
\def\e{\eta}
\def\ep{\varepsilon}
\def\G{\Gamma}
\def\g{\gamma}

\def\la{\lambda}
\def\La{\Lambda}
\def\om{\omega}
\def\Om{\Omega}
\def\th{\theta}

\def\w{\wedge}

\begin{document}
\title
[curvature estimates  for minimal submanifolds] {curvature
estimates  for minimal submanifolds  of higher codimension}

\author
[Y.L. Xin and Ling Yang]{Y. L. Xin and Ling Yang}
\address
{Institute of Mathematics, Fudan University, Shanghai 200433, China
and Key Laboratory of Mathematics for Nonlinear Sciences (Fudan
University), Ministry of Education} \email{ylxin@fudan.edu.cn}
\thanks{The research was partially supported by
NSFC}
\begin{abstract}
We derive curvature estimates for minimal submanifolds in Euclidean
space for arbitrary dimension and codimension via Gauss map. Thus,
Schoen-Simon-Yau's results and Ecker-Huisken's results are
generalized to higher codimension. In this way we improve
Hildebrandt-Jost-Widman's result for the Bernstein type theorem.
\end{abstract}

\renewcommand{\subjclassname}{%
  \textup{2000} Mathematics Subject Classification}
\date{}
\maketitle

\Section{Introduction}{Introduction}
\medskip
Let $f$ be a smooth function on an open domain $\Om\subset\ir{n}.$
 If $f$ satisfies the minimal surface equation, it defines a minimal
hypersurface $M$ in $\ir{n+1}$.If $f$ is an entire solution to the
equation, $f$ must be an affine linear function for $n\le 7$ whose
graph is a hyperplane . Those are the classical Bernstein theorem
\cite{b} and its higher dimensional generalizations which was
finally proved by J. Simons \cite{si} . Counterexamples to the
theorems for $n\ge 8$ were given by Bombeiri-De Gorge-Guisti
\cite{b-d-g}.

Heinze \cite{h} considered the minimal graph defined over a disc
$D_R\subset \ir{2}$ and gave  curvature estimates. The classical
Benstein theorem can be obtained by letting $R\to +\infty$ in his
curvature estimates.

For general minimal surface in Euclidean space, so-called parametric
case, the Beinstein type results are closely related to the value
distribution of the Gauss image. The work of
Osserman-Xavier-Fujimoto \cite{o}\cite{xa}\cite{f} has  settled the
question of the value of Gauss map for complete minimal surfaces in
$\ir{3}$.

A minimal graph is area-minimizing, in particular, it is stable in
the sense that its second variation of the volume is non-negative on
any compact subset of $M$. Any stable minimal surface in $\ir{3}$ is
a plane, due to Fisher- Colbrie- Schoen and do Carmo-Peng
\cite{f-c-s}\cite{d-p}. But, for higher dimensional stable
hypersurfaces  Cao-Shen-Zhu proved that they have only one end
\cite{c-s-z}.

For stable minimal hypersurfaces, Schoen-Simon-Yau gave curvature
estimates,  which not only gave us a direct proof for Bernstein type
theorems for $n\le 5$ dimensional minimal graphs, but also gave us a
new method to obtain curvature estimates.

For any $n\ge 2$, there is a weak version of the Bernstein type
theorem. It was J. Moser \cite{m} who proved that the entire
solution $f$ to the minimal surface equation is affine linear,
provided $|\n f|$ is uniformly bounded. There is no dimension
limitation. By a geometric approach Ecker-Huisken \cite{e-h}
obtained the curvature estimates, as a conclusion Moser's result has
been improved  for the controlled growth of $|\n f|$.

For area-minimizing hypersurfaces  with vanishing first Betti number
, Solomom \cite{so} was able to give curvature estimates under the
hypotheses of Gauss map.

Higher codimensional Bernstein problem becomes more complicated.
There is a counterexample given by Lawson-Osserman \cite{l-o}. On
the other hand,  Hildebrandt-Jost-Widman \cite{h-j-w} generalized
Moser's result to higher codimension as follows.

\begin{thm}\label{C} Let $z^\a=f^\a(x), \a=1,\cdots,m,\;x=(x^1,\cdots,x^n)\in
\ir{n}$ be the $C^2$ solution to the system of minimal surface
equations. Let there exist $\be,$ where
$$\be<\cos^{-p}\left(\frac \pi{2\sqrt{p\,K}}\right),\quad K=\begin{cases} 1 \quad\text{if}\quad p=1\cr
         2 \quad \text{if} \quad p\ge 2\end{cases},\quad p=\min(m, n)$$
such that for any $x\in \ir{n}$,
$$\De_f(x)
=\{\text{det}(\de_{ij}+f_{x^i}^s(x)f_{x^j}^s(x))\}^{\frac
12}\le\be,$$ then $f^1,\cdots,f^m$ are affine linear functions on
$\ir{n},$ whose graph is affine $n-$plane in $\ir{m+n}.$
\end{thm}

The geometric meaning of the condition  in the above theorem is that
the image under the Gauss map lies in a closed subset of an open
geodesic ball of the radius $\f{\sqrt{2}}{4}\pi$. Later in a joint
work of the first author with J. Jost \cite{j-x} we found larger
geodesic convex set $B_{JX}(P_0)$, where $P_0$ denotes a fixed
$n-$plane, and then improved the above theorem. Our bound of slope
is $2$, big than $\cos^{-p}\left(\frac \pi{2\sqrt{p\,K}}\right)$. It
should be noted that although $B_{JX}(P_0)\supset
B_{\f{\sqrt{2}}{4}\pi}(P_0)$,  they have some common boundary
points.

Recently, the first author and his collaborators
\cite{x3}\cite{s-w-x} studied complete minimal submanifolds with
flat normal bundle and positive $w-$function. In this special
situation, the Schoen-Simon-Yau type curvature estimates and the
Ecker-Huisken type curvature estimates can be carried out, then the
corresponding Bernstein type theorems follow immediately.

In this paper we  study a complete minimal submanifold $M$ in
$\ir{m+n}$ with the codimension $m\ge 2$. We have a Bochner type
formula for the squared norm of the second fundamental form $B$. As
in the codimension one case, we need a Kato type inequality for $|\n
B|^2$ in terms of $|\n |B||^2$. This is one of the crucial issue in
the curvature estimates. We derive it for any codimension in \S 2.

Since the curved normal bundle, the curvature estimates would be
more delicate. We can define Gauss map from $M$ to the Grassmannian
manifold $\grs{n}{m}$. From the counterexample of Lawson-Osserman,
some additional conditions are needed to study higher codimensional
Bernstein problem. The adequate conditions would  confine the image
of the Gauss map, as in the previous work of
Osserman-Xavier-Fujimoto in dimension $2$ and in the work of Solomon
in higher dimension. Now, in general dimension and codimension, we
assume that the image under the Gauss map lies in an open geodesic
ball of radius $\f{\sqrt{2}}{4}\pi$ in $\grs{n}{m}$ which is the
largest convex geodesic ball in Grassmannian. We find two auxiliary
functions $h_1$ and $h_2$ on $M$ via Gauss map. The precise
definition and their properties can be found in \S 3.

(\ref{dh1}) shows that $h_1$ can be viewed as a generalized 
support function.  With the aid of the function $h_1$, we can 
derive a "strong stability inequality" which enables us to carry 
out the curvature estimates of Schoen-Simon-Yau type.

Using the function $h_2$ we can find subharmonic functions on $M$
and carry out the curvature estimates of Ecker-Huisken type in terms
of $h_2$.

By those curvature estimates we can obtain the following Bernstein
type theorems.
\begin{thm}\label{A}
Let $M$ be a complete minimal $n$-submanifold in $\ir{n+m}$ ($n\leq
6$). If the Gauss image of $M$ is contained in an open geodesic ball
of $\grs{n}{m}$ centered at $P_0$ and of radius
$\f{\sqrt{2}}{4}\pi$, then $M$ has to be an affine linear subspace.

\end{thm}

\begin{thm}\label{B}
Let $M$ be a complete minimal $n$-submanifold in $\ir{n+m}$. If the
Gauss image of $M$ is contained in an open geodesic ball of
$\grs{n}{m}$ centered at $P_0$ and of radius $\f{\sqrt{2}}{4}\pi$,
and $\big(\f{\sqrt{2}}{4}\pi-\rho\circ \g\big)^{-1}$ has growth
\begin{equation}
\big(\f{\sqrt{2}}{4}\pi-\rho\circ \g\big)^{-1}=o(R),
\end{equation}
where $\rho$ denotes the distance on $\grs{n}{m}$ from $P_0$ and $R$ is
the Euclidean distance from any point in $M$. Then $M$ has to be an
affine linear subspace.
\end{thm}

In general, if there is neither dimension limitation as in Theorem
\ref{A} nor the growth assumption as in Theorem \ref{B}, we conclude
that $M$ has only one end which is a consequence of the "strong
stability inequality".

Our Theorem \ref{A} and Theorem \ref{B} are closely related to
Theorem \ref{C} mentioned above.  Our results are higher
codimensional generalizations of Schoen-Simon-Yau's results and
Ecker-Huisken's results, and improve Hildebrandt-Jost-Widman's
theorem.

It is worthy to note that our method is also suitable for
codimension one case. We only need to modify the auxiliary functions
$h_1$ and $h_2$ in order to recover the known results for minimal
hypersurfaces.
\bigskip\bigskip

\Section{A Kato-Type Formula}{A Kato-Type Formula}
\medskip

Let $M\to \bar M$ be an isometric immersion with the second
fundamental form $B,$ which can be viewed as a cross-section of the
vector bundle Hom($\odot^2TM, NM$) over $M,$ where  $TM$ and $NM$
denote the tangent bundle and the normal bundle  along $M$,
respectively. A connection on Hom($\odot^2TM, NM$) can be induced
from those of $TM$ and $NM$ naturally.

For $\nu\in\G(NM)$ the shape operator $A^\nu: TM\to TM$  satisfies
$$\left<B_{X Y}, \nu\right>=\left<A^\nu(X), Y\right>.$$

We define the mean curvature $H$ to be the trace of the second
fundamental form. It is a normal vector field on $M$ in $\bar M$.

The second fundamental form, curvature tensor of the submanifold,
curvature tensor of the normal bundle and that of the ambient
manifold satisfy the Gauss equations, the Codazzi equations and the
Ricci equations as follows.
$$\left<R_{X Y}Z, W\right>
=\left<\bar R_{X Y}Z, W\right>-\left<B_{X W}, B_{Y Z}\right>
                                   +\left<B_{X Z}, B_{Y W}\right>,$$
$$(\n_XB)_{Y Z}-(\n_YB)_{X Z}=-(\bar R_{X Y}Z)^N,$$
$$\left<R_{X Y}\mu, \nu\right>=\left<\bar R_{X Y}\mu, \nu\right>
  +\left<B_{X e_i}, \mu\right>\left<B_{Y e_i}, \nu\right>
    -\left<B_{X e_i}, \nu\right>\left<B_{Y e_i}, \mu\right>,$$
where $\{e_i\}$ is a local orthonormal frame field of $M$;
$X,\,Y\,$ and $Z$ are tangent vector fields; $\mu,\, \nu$ are
normal vector fields in $M.$ Here and in the sequel we use the
summation convention and  agree the range of indices:
$$1\leq i, j, s, t\leq n;\qquad 1\leq \a, \be, \g\leq m;\qquad 1\le a, b, c\le m+n.$$

There is the trace-Laplace operator $\n^2$ acting on any
cross-section of a Riemannian vector bundle $E$ over $M$.

Now, we consider a minimal submanifold $M$ of dimension $n$ in
Euclidean ($m+n$)-space $\ir{m+n}$ with $m\ge 2$.  We have (see
\cite{si})
\begin{equation}
\n^2B=-\tilde{\mc{B}}-\underline{\mc{B}}.\label{Si} \end{equation}
We recall the following notations:
$$\tilde{\mc{B}}\mathop{=}\limits^{def.}B\circ B^t\circ B,$$
where $B^t$ is the conjugate map of $B;$
$$\underline{\mc{B}}_{X\, Y}\mathop{=}\limits^{def.}
\sum_{\a=1}^m\left(B_{A^{\nu_\a}A^{\nu_\a}(X)\,Y}+B_{X\,A^{\nu_\a}A^{\nu_\a}(Y)}
-2\,B_{A^{\nu_\a}(X)\,A^{\nu_\a}(Y)}\right),$$ where $\nu_\a$ are
basis vectors of normal space. It is obvious that
$\underline{\mc{B}}_{X\, Y}$ is symmetric in $X$ and $Y,$ which is a
cross-section of the bundle Hom($\odot^2TM, NM$). Simons also gave
an estimate \cite{si}
$$\left<\tilde{\mc{B}}+\underline{\mc{B}}, B\right>
\le \left(2-\frac 1m\right)|B|^4.$$ It is optimal for the
codimension $m=1.$

In the case when $m\ge 2,$ there is a refined estimate
\cite{l-l}\cite{c-x}
$$\left<\tilde{\mc{B}}+\underline{\mc{B}}, B\right> \le \frac 32 |B|^4.$$
Substituting it into (\ref{Si}) gives
$$\left<\n^2B,B\right>\ge -\frac 32|B|^4.$$
It follows that
\begin{equation}\De|B|^2\ge -3 |B|^4+2|\n B|^2.\label{bo}\end{equation}

We need a Kato-type inequality in order to use the  formula
(\ref{bo}). Namely, we would estimate $|\n B|^2$ in terms of
$|\n|B||^2.$ Schoen-Simon-Yau \cite{s-s-y} did such an estimate for
codimension $m=1.$ For any $m$ with flat normal bundle their
technique is also applicable, see the previous paper of the first
author \cite{x3}. The following lemma is a generalized version of
their estimate for any codimension $m$.

\begin{lem}
\begin{equation}
|\n B|^2\geq \big(1+\f{2}{mn}\big)\big|\n|B|\big|^2\label{kato}.
\end{equation}
\end{lem}
\begin{proof} It is sufficient for us to prove the inequality at the points
where $|B|^2\neq 0$. Choose a local orthonormal tangent frame field
$\{e_1,\cdots, e_n\}$ and a local orthonormal normal frame field
$\{\nu_1,\cdots,\nu_m\}$ of $M$ near the considered point $x$, such
that
\begin{equation}\label{co15}
\n_{e_i}e_j(x)=0,\ \n_{e_i}\nu_\a(x)=0,
\end{equation}
Denote the shape operator $A^\a=A^{\nu_\a}.$  Then obviously
$|B|^2=\sum_\a |A^\a|^2$ and
\begin{eqnarray*}
\n |B|^2&=&\sum_\a \n |A^\a|^2.\\
\big|\n |B|^2\big|^2&=&\sum_{\a,\be} \n |A^\a|^2\cdot \n |A^\be|^2
\leq \f{1}{2}\sum_{\a,\be}\Big(\big|\n |A^\a|^2\big|^2+\big|\n |A^\be|^2\big|^2\Big)\\
&\leq&m\sum_\a \big|\n |A^\a|^2\big|^2.
\end{eqnarray*}
Therefore
\begin{equation}\label{ineq7}
\big|\n|B|\big|^2=\f{\big|\n |B|^2\big|^2}{4|B|^2}\leq \f{m\sum_\a
\big|\n |A^\a|^2\big|^2}{4\sum_\a |A^\a|^2}.
\end{equation}

Since $|B|^2\neq 0$, there exist $\g,k,l$ such that $h_{\g kl}\neq
0$, where $h_{\a ij}=\left< B_{e_ie_j},\nu_\a\right>$ for arbitrary
$\a,\  i,\ j$; then
$$(h_{1kl},\cdots,h_{mkl})\in \ir{m}-\{0\}.$$
Obviously, there exists an $m\times m$ orthogonal matrix $U$ and
$z=(z_1,\cdots,z_m)\in \ir{m}$, such that
\begin{equation*}
z_\a\neq 0\qquad \mbox{for every }1\leq \a\leq m
\end{equation*}
and
\begin{equation*}
z_\a=U_\a^\be h_{\be kl}.
\end{equation*}
Now we define $\td{\nu}_\a=U_\a^\be \nu_\be$, then
$\{\td{\nu}_1,\cdots,\td{\nu}_m\}$ is also a local orthogonal normal
frame field, satisfying (\ref{co15}); and moreover
\begin{equation*}
\td{h}_{\a kl}=\left< B_{e_ke_l},\td{\nu}_\a\right>=U_\a^\be h_{\be
kl}=z_\a\neq 0.
\end{equation*}
Define the shape operator $\td{A}^\a$ corresponding to $\td\nu_\a$,
then
\begin{equation*}
|\td{A}^\a|^2=\sum_{i,j}\td{h}_{\a ij}^2\geq \td{h}_{\a kl}^2>0.
\end{equation*}
Hence we can assume $|A^\a|^2>0$ for arbitrary $\a$ without loss of
generality.

Let $1\leq \g\leq m$ such that
\begin{equation*}
\f{\big|\n |A^\g|^2\big|^2}{|A^\g|^2}=\max_\a \Bigg\{\f{\big|\n
|A^\a|^2\big|^2}{|A^\a|^2}\Bigg\}<+\infty,
\end{equation*}
then from (\ref{ineq7}),
\begin{equation}
\big|\n|B|\big|^2\leq \f{m \big|\n |A^\g|^2\big|^2}{4|A^\g|^2}.
\end{equation}
$|A^\g|^2$ and $\n|A^\g|^2$ is independent of the choice of
$\{e_1,\cdots,e_n\}$, then without loss of generality we can assume
$h_{\g ij}=0$ whenever $i\neq j$. Then
\begin{eqnarray*}
\n|A^\g|^2&=&2\sum_k\sum_{i,j}h_{\g ij}h_{\g ijk}e_k=2\sum_k\sum_i h_{\g ii}h_{\g iik}e_k,\nonumber\\
\big|\n|A^\g|^2\big|^2
&=&4\sum_k \big(\sum_i h_{\g ii}h_{\g iik}\big)^2
\leq 4\sum_k \big(\sum_i h_{\g ii}^2\big)\big(\sum_i h_{\g iik}^2\big)\nonumber\\
&=&4\big(\sum_i h_{\g ii}^2\big)\big(\sum_{i,k}h_{\g
iik}^2\big)=4|A^\g|^2\sum_{i,k}h_{\g iik}^2\nonumber
\end{eqnarray*}
and
\begin{eqnarray}
\big|\n |B|\big|^2&\leq& \f{m \big|\n |A^\g|^2\big|^2}{4|A^\g|^2}\leq m\sum_{i,k}h_{\g iik}^2\nonumber\\
&=&m\sum_{i\neq k} h_{\g iik}^2+m\sum_i h_{\g iii}^2\nonumber\\
&=&m\sum_{i\neq k} h_{\g iki}^2+m\sum_i (\sum_{j\neq i}h_{\g jji})^2\nonumber\\
&\leq&m\sum_{i\neq k}h_{\g iki}^2+(n-1)m\sum_{i\neq j}h_{\g jji}^2\nonumber\\
&=&nm\sum_{i\neq k}h_{\g iki}^2,\label{ineq8}
\end{eqnarray}
where we used the Codazzi equations and the vanishing mean
curvature condition $H=0$. Please note that $h_{\a ijk}=\left<
(\n_{e_k} B)_{e_i,e_j},\nu_\a\right>$ for arbitrary $\a, i, j, k.$

On the other hand, a direct calculation shows
\begin{eqnarray}
\big|\n|B|^2\big|^2&=&|2\sum_k \sum_{\a,i,j}h_{\a ij}h_{\a ijk}e_k|^2
=4\sum_{\a,\be,i,j,s,t,k}h_{\a ij}h_{\a ijk}h_{\be st}h_{\be stk},\nonumber\\
|\n B|^2-\big|\n|B|\big|^2&=&|\n B|^2-\f{\big|\n|B|^2\big|^2}{4|B|^2}\nonumber\\
&=&\sum_{\a,i,j,k}h_{\a ijk}^2
-\f{\sum_{\a,\be,i,j,s,t,k}h_{\a ij}h_{\a ijk}h_{\be st}h_{\be stk}}{\sum_{\be,s,t} h_{\be st}^2}\nonumber\\
&=&\f{\sum_{\a,i,j,s,t,k}(h_{\a ijk}h_{\be st}-h_{\be stk}h_{\a ij})^2}{2|B|^2}\nonumber\\
&\geq&\f{\sum_{\be,i\neq j,s,t,k}h_{\g ijk}^2h_{\be st}^2
+\sum_{\a,s\neq t,i,j,k}h_{\g stk}^2 h_{\a ij}^2}{2|B|^2}\nonumber\\
&=&\sum_{i\neq j,k}h_{\g ijk}^2\geq \sum_{i\neq k}(h_{\g iki}^2+h_{\g ikk}^2)\nonumber\\
&=&2\sum_{i\neq k}h_{\g iki}^2.
\end{eqnarray}
In conjunction with (\ref{ineq8}), we finally arrive at
(\ref{kato}).
\end{proof}

Hence it follows  from (\ref{bo}) and  (\ref{kato}) that
\begin{equation}\label{ineq9}
\De |B|^2\geq 2\big(1+\f{2}{mn}\big)\big|\n |B|\big|^2-3|B|^4.
\end{equation}
\bigskip\bigskip

\Section{Auxiliary Functions via Gauss Maps}{Auxiliary Functions via
Gauss Maps}
\medskip

Let $\ir{m+n}$ be an $(m+n)$-dimensional Euclidean space. All
oriented $n$-subspaces constitute the Grassmannian manifold
$\grs{n}{m}$, which is an irreducible symmetric space of compact
type.  The canonical Riemannian metric on $\grs{n}{m}$ can be
expressed in the following way.

Let $\{e_i,e_{n+\a}\}$ be a local orthonormal frame field in
$\ir{m+n}$. Let $\{\om_i,\om_{n+\a}\}$ be the dual frame field so
that the Euclidean metric is
\begin{equation*}
g=\sum_i \om_i^2+\sum_\a \om_{n+\a}^2.
\end{equation*}
The Levi-Civita connection forms $\om_{ab}$ of $\ir{m+n}$ are
uniquely determined by the equations
\begin{eqnarray*}
&&d\om_a=\om_{ab}\w \om_b,\nonumber\\
&&\om_{ab}+\om_{ba}=0.
\end{eqnarray*}
Let $P\in \grs{n}{m}$ be any point which is spanned by
$\{e_1,\cdots,e_n\}$, then the canonical Riemannian metric on
$\grs{n}{m}$ can be written as
\begin{equation}\label{a}
ds^2=\sum_{i,\a}\om_{n+\a\ i}^2.
\end{equation}
The sectional curvature of the above canonical metric varies in
the interval $[0, 2]$ in the case of $\min\{ n, m\}\ge 2$. By the
standard Hessian comparison theorem we have
\begin{equation}
\text{Hess}(\rho)\ge \sqrt{2}\,\cot(\sqrt{2}\rho)(g-d \rho\otimes d
\rho),\label{hr}
\end{equation}
where $\rho$ is the distance function from a fixed point in
$\grs{n}{m}$ and $g$ is the metric tensor on $\grs{n}{m}.$

Let $0$ be the origin of  $\ir{m+n}$. Let $ SO (m+n)$ be the Lie
group consisting of all the orthonormal frames $ (0; e_i, e_{n+\a})
$. Let $ P = \{ (x; e_1, ..., e_n ): x \in M, e_i \in T_x M \} $ be
the principal bundle of orthonormal tangent frames over $ M, Q =\{ (
x; e_{n+1}, ..., e_{n+m}  ) : x \in M, e_{n+\a} \in N_x M \} $ be
the principal bundle of orthonormal normal frames over $ M$, then $
\bar\pi : P \oplus Q \to M $ is the projection with fiber $ SO (m)
\times SO (n), \;i: P \oplus Q \hookrightarrow SO (m+n)  $ is the
natural inclusion.

The  Gauss map $ \gamma : M \to \grs{n}{m} $ is defined by
$$
 \g (x) = T_x M \in \grs{n}{m}
$$
via the parallel translation in $ \ir{m+n} $ for $ \forall x \in M
$. Thus, the following commutative  diagram holds

$$\CD
 P \oplus Q @>i>> SO (m+n)  \\
 @V{\bar\pi}VV     @VV{\pi}V \\
 M  @>{\g}>>  \grs{n}{m}
\endCD$$

From the above diagram we know that the energy density of the Gauss
map (see \cite{x1} Chap.3, \S 3.1)
$$e(\g)=\f{1}{2}\left<\g_*e_i,\g_*e_i\right>=\f{1}{2}|B|^2.$$
Ruh-Vilms proved  that the mean curvature vector of $M$ is parallel
if and only if its Gauss map is a harmonic map \cite{r-v}.

We consider smooth functions on an open geodesic ball
$B_{\f{\sqrt{2}}{4}\pi}(P_0)\subset\grs{n}{m}$ of radius
$\frac{\sqrt{2}}{4}\pi$ and centered at $P_0$. Those are useful for
our curvature estimates later.  Let
$$u=\cos(\sqrt{2}\rho),$$
where $\rho$  is the distance function from $P_0$ in $\grs{n}{m}$.
We have
$$u'=-\sqrt{2}\sin(\sqrt{2}\rho),$$
$$u''=-2\cos(\sqrt{2}\rho).$$
Then,
\begin{equation}\aligned
\text{Hess}(u)&=u'\text{Hess}(\rho)+u''d\rho\otimes d\rho\\
&\le -2\cos(\sqrt{2}\rho)(g-d\rho\otimes
d\rho)-2\cos(\sqrt{2}\rho)d\rho\otimes d\rho=-2 u g.
\endaligned
\end{equation}

The composition function $h_1=u\circ \g$ of $u$ with the Gauss map
$\g$ defines a function on $M$. Using the composition formula, we
have
\begin{equation}
\aligned \De h_1&=\text{Hess}(u)(\g_* e_i, \g_* e_i)+d u(\tau(\g))\\
&\le -2 |B|^2h_1,
\endaligned\label{dh1}
\end{equation}
where $\tau(\g)$ is the tension field of the Gauss map, which is
zero, provided $M$ has parallel mean curvature by the Ruh-Vilms
theorem mentioned above.

Let
$$h=\sec^2(\sqrt{2}\rho),$$
where $\rho$ is the distance function from $P_0$ in $\grs{n}{m}$. We
have
$$h'=2 \sqrt{2}\sec^2(\sqrt{2}\rho)\tan(\sqrt{2}\rho),$$
$$h''=12\ \sec^2(\sqrt{2}\rho)\tan^2(\sqrt{2}\rho)+4\ \sec^2(\sqrt{2}\rho).$$
Hence,
$$\aligned
\text{Hess}(h)&=h'\text{Hess}(\rho)+h''d\rho\otimes d\rho\\
&\ge 4\ \sec^2(\sqrt{2}\rho)(g-d\rho\otimes
d\rho)\\
&\qquad+\left( 12\sec^2(\sqrt{2}\rho)\tan^2(\sqrt{2}\rho)+ 4
\sec^2(\sqrt{2}\rho)\right)d\rho\otimes
d\rho\\
&= 4 h\ g+\f{3}{2}h^{-1}dh\otimes dh\endaligned $$

The composition function $h_2=h\circ \g$ of $h$ with the Gauss map
$\g$ defines a function on $M$. Using the composition formula, we
have
\begin{equation}
\aligned \De h_2&=\text{Hess}(h)(\g_* e_i, \g_* e_i)+d h(\tau(\g))\\
&\ge 4\ h_2 |B|^2+\f{3}{2}h_2^{-1}|\n h_2|^2,
\endaligned\label{dh2}
\end{equation}
where $\tau(\g)$ is the tension field of the Gauss map, which is
zero in our consideration.

With the aid of $h_1$, we immediately have the following lemma.
\begin{lem}
Let $M$ be an $n$-dimensional minimal submanifold of $\ir{n+m}$
($M$ needs not be complete), if the Gauss image of $M$ is
contained in an open geodesic ball of radius $\f{\sqrt{2}}{4}\pi$
in $\grs{n}{m}$, then we have
\begin{equation}\label{lp3}
\int_M |\n \phi|^2*1\geq 2\int_M |B|^2\phi^2*1
\end{equation}
for any function $\phi$ with compact support $D\subset M$.
\end{lem}
\begin{proof}
Let
$$L\phi=-\De \phi-2|B|^2\phi.$$
Its first eigenvalue with the Dirichlet boundary condition in $D$ is
$\la_1$ and the corresponding eigenfunction is $v$. Without loss of
generality, we assume that $v$ achieves the positive maximum.
Consider a $C^2$ function
$$f=\f{v}{h_1}.$$
Since $f|_{\p D}=0$, it achieves the positive maximum at a point
$x\in D$. Therefore, at $x$,
$$\n f=0,\qquad \De f\leq 0.$$
It follows that
\begin{eqnarray*}
\De v&=&\De(f h_1)=\De f\cdot h_1 + f\De h_1+2\n f\cdot\n h_1\\
&\leq&f\De h_1=\f{v\De h_1}{h_1}.
\end{eqnarray*}
Namely, at $x$,
$$\f{\De v}{v}\leq \f{\De h_1}{h_1}.$$
\begin{equation}\label{lp1}
\f{\De v+2|B|^2v}{v}\leq \f{\De h_1+2|B|^2 h_1}{h_1}\leq 0.
\end{equation}
On the other hand,
\begin{equation}\label{lp2}
\f{\De v+2|B|^2v}{v}=-\la_1.
\end{equation}
(\ref{lp1}) and (\ref{lp2}) implies $\la_1\geq 0$. Hence we have
$$0\leq \la_1=\inf\f{\int_D \phi L\phi*1}{\int_D \phi^2*1}\leq \f{\int_D \phi L\phi*1}{\int_D \phi^2*1},$$
which shows (\ref{lp3}) holds true.
\end{proof}

\begin{rem}
For a stable minimal hypersurface there is the stability inequality,
which is one of main ingredient for Schoen-Simon-Yau's curvature
esimates for stable minimal hypersurfaces. For minimal submanifolds
with the Gauss image restriction we have stronger inequality as
shown in (\ref{lp3})
\end{rem}
\bigskip\bigskip

\Section{Curvature estimates of Schoen-Simon-Yau type}{Curvature
estimates of Schoen-Simon-Yau type}
\medskip

Replacing $\phi$ by $|B|^{1+q}\phi$ in (\ref{lp3}) gives
\begin{equation}\label{lp}\aligned
\int_M |B|^{4+2q}&\phi^2*1\leq \f{1}{2}\int_M \big|\n (|B|^{1+q}\phi)\big|^2*1\\
&=\f{1}{2}(1+q)^2\int_M |B|^{2q}\big|\n|B|\big|^2\phi^2*1+\f{1}{2}\int_M |B|^{2+2q}|\n \phi|^2*1\\
&+(1+q)\int_M |B|^{1+2q}\n|B|\cdot \phi\n\phi *1.\endaligned
\end{equation}
From (\ref{ineq9})  we can derive
\begin{equation}\label{lp5}
\f{2}{mn}\big|\n|B|\big|^2\leq |B|\De |B|+\f{3}{2}|B|^4.
\end{equation}
Multiplying $|B|^{2q}\phi^2$ with both sides of (\ref{lp5}) and
integrating by parts, we have
\begin{equation}\label{lp6}\aligned
\f{2}{mn}\int_M |B|^{2q}&\big|\n |B|\big|^2\phi^2*1
\leq\int_M |B|^{1+2q}\De |B|\phi^2*1+\f{3}{2}\int_M |B|^{4+2q}\phi^2*1\\
&=-\int_M \n|B|\cdot\n\big(|B|^{1+2q}\phi^2\big)*1+\f{3}{2}\int_M |B|^{4+2q}\phi^2 *1\\
&=-(1+2q)\int_M |B|^{2q}\big|\n|B|\big|^2\phi^2*1\\
&\qquad -2\int_M |B|^{1+2q}\n|B|\cdot \phi\n\phi*1 +\f{3}{2}\int_M
|B|^{4+2q}\phi^2 *1. \endaligned
\end{equation}

By multiplying $\f{3}{2}$ with both sides of (\ref{lp}) and then
adding up both sides of it and (\ref{lp6}), we have
\begin{equation}\label{lp7}\aligned
&\big(\f{2}{mn}+1+2q-\f{3}{4}(1+q)^2\big)\int_M |B|^{2q}\big|\n |B|\big|^2\phi^2*1\\
&\leq\f{3}{4}\int_M |B|^{2+2q}|\n\phi|^2*1
+\big(\f{3}{2}(1+q)-2\big)\int_M |B|^{1+2q}\n|B|\cdot \phi\n\phi
*1.
\endaligned\end{equation} By using Young's inequality,
\begin{equation*}\aligned
&\Big(\f{3}{2}(1+q)-2\Big)\int_M |B|^{1+2q}\n|B|\cdot \phi\n\phi *1\\
&\qquad\leq \ep \int_M |B|^{2q}\big|\n
|B|\big|^2\phi^2*1+C_1(\ep,q)\int_M |B|^{2+2q}|\n\phi|^2*1.
\endaligned\end{equation*} Then (\ref{lp7}) becomes
\begin{equation}\label{Lp8}\aligned
&\big(\f{2}{mn}+1+2q-\f{3}{4}(1+q)^2-\ep\big)\int_M |B|^{2q}\big|\n |B|\big|^2\phi^2*1\\
&\hskip1in\leq C_2(\ep,q)\int_M |B|^{2+2q}|\n\phi|^2*1.\endaligned
\end{equation}
When
\begin{equation}\label{Lp14} q\in
\Big[0,\f{1}{3}+\f{2}{3}\sqrt{1+\f{6}{mn}}\ \Big)
\end{equation}
we have
$$\f{2}{mn}+1+2q-\f{3}{4}(1+q)^2>0.$$
Then we can choose $\ep$ sufficiently small, such that
\begin{equation}\label{Lp9}
\int_M |B|^{2q}\big|\n |B|\big|^2\phi^2*1\leq C_3 \int_M
|B|^{2+2q}|\n\phi|^2*1
\end{equation}
where $C_3$ only depends on $n$, $m$ and $q$.

Again using Young's inequality yields
\begin{equation}\label{lp10}
|B|^{1+2q}\n|B|\cdot \phi\n\phi\leq
\f{1}{2}\big(|B|^{2q}\big|\n|B|\big|^2\phi^2+|B|^{2+2q}|\n\phi|^2\big).
\end{equation}
Substituting (\ref{Lp9}) and (\ref{lp10}) into (\ref{lp}) gives
\begin{equation}\label{Lp11}
\int_M |B|^{4+2q}\phi^2*1\leq C_4(n,m,q)\int_M |B|^{2+2q}|\n
\phi|^2*1.
\end{equation}

Replacing $\phi$ by $\phi^{2+q}$ in (\ref{Lp11}) gives
$$\int_M |B|^{4+2q}\phi^{4+2q}*1\leq C_4(2+q)^2\int_M |B|^{2+2q}\phi^{2+2q}|\n \phi|^2*1.$$
By H\"{o}lder inequality, we have
\begin{eqnarray*}
&&\int_M |B|^{2+2q}\phi^{2+2q}|\n \phi|^2*1\\
&\leq&\big(\int_M
|B|^{4+2q}\phi^{4+2q}*1\big)^{\f{1+q}{2+q}}\big(\int_M
|\n\phi|^{4+2q}*1\big)^{\f{1}{2+q}}.
\end{eqnarray*}
Therefore
\begin{eqnarray}\label{Lp12}
\int_M |B|^{4+2q}\phi^{4+2q}*1\leq C\int_M |\n\phi|^{4+2q}*1.
\end{eqnarray}
where $C$ is a constant only depending on $n$, $m$ and $q$.

Replacing $\phi$ by $\phi^{1+q}$ in (\ref{Lp11}) gives
$$\int_M |B|^{4+2q}\phi^{2+2q}*1\leq C_4(1+q)^2\int_M |B|^{2+2q}\phi^{2q}|\n \phi|^2*1.$$
By H\"{o}lder inequality, we have
\begin{eqnarray*}
&&\int_M |B|^{2+2q}\phi^{2q}|\n \phi|^2*1\\
&\leq&\big(\int_M
|B|^{4+2q}\phi^{2+2q}*1\big)^{\f{q}{1+q}}\big(\int_M
|B|^2|\n\phi|^{2+2q}*1\big)^{\f{1}{1+q}}.
\end{eqnarray*}
Therefore
\begin{eqnarray}\label{Lp13}
\int_M |B|^{4+2q}\phi^{2+2q}*1\leq C'\int_M |B|^2|\n\phi|^{2+2q}*1.
\end{eqnarray}
where $C'$ is a constant only depending on $n$, $m$ and $q$.

Let $r$ be a function on $M$ with $|\n r|\le 1$. For any $R\in [0,
R_0]$, where $R_0=\sup_Mr$, suppose
$$M_R=\{x\in M,\quad r\le R\}$$
is compact.

(\ref{Lp12}) enables us to prove the following results.

\begin{thm}\label{t1}
Let $M$ be an $n$-dimensional minimal submanifolds of $\ir{n+m}$.
If  the Gauss image of $M_R$ is contained in an open geodesic
ball of radius $\f{\sqrt{2}}{4}\pi$ in $\grs{n}{m}$, then we have
the $L^p$-estimate
\begin{equation}\label{Lp15}
\big\||B|\big\|_{L^p(M_{\th R})}\leq
C(n,m,p)(1-\th)^{-1}R^{-1}\text{Vol}(M_R)^{\f{1}{p}}
\end{equation}
for arbitrary $\th\in (0,1)$ and
$$p\in \left[4,4+\f{2}{3}+\f{4}{3}\sqrt{1+\f{6}{mn}} \right).$$
\end{thm}

\begin{proof} Take $\phi\in C_c^\infty (M_R)$ to be the
standard cut-off function such that $\phi\equiv 1$ in $M_{\th R}$
and  $|\n \phi|\leq C(1-\th)^{-1}R^{-1}$; then (\ref{Lp12}) yields
$$\int_{M_{\th R}}|B|^p*1\leq C(1-\th)^{-p}R^{-p}\text{Vol}(M_R),$$
where $p=4+2q$. Thus the conclusion immediately follows from
(\ref{Lp12}).
\end{proof}

\Section{Curvature estimates of Ecker-Huisken type}{Curvature
estimates of Ecker-Huisken type}
\medskip

We compute from (\ref{ineq9}) and (\ref{dh2}):
\begin{eqnarray*}
&&\De\big(|B|^{2p}h_2^q\big)\\
&=&\De |B|^{2p}\cdot h_2^q+|B|^{2p}\De h_2^q+2\n |B|^{2p}\cdot \n h_2^q\\
&=&\Big(p|B|^{2p-2}\De |B|^2+p(p-1)|B|^{2p-4}\big|\n|B|^2\big|^2\Big)h_2^q\\
&&+|B|^{2p}\big(q h_2^{q-1}\De h_2+q(q-1)h_2^{q-2}|\n h_2|^2\big)+4pq|B|^{2p-1}\n|B|\cdot h_2^{q-1}\n h_2\\
&\geq& (4q-3p)|B|^{2p+2}h_2^q+2p(2p-1+\f{2}{mn})|B|^{2p-2}\big|\n|B|\big|^2 h_2^q
+q(q+\f{1}{2})|B|^{2p}h_2^{q-2}|\n h_2|^2\\
&&+4pq|B|^{2p-1}\n|B|\cdot h_2^{q-1}\n h_2.
\end{eqnarray*}
By Young's inequality, when $2p(2p-1+\f{2}{mn})\cdot
q(q+\f{1}{2})\geq (2pq)^2$, i.e.,
\begin{equation}
p\geq \f{1}{2}-\f{1}{mn}+\big(1-\f{2}{mn}\big)q,\label{p1}
\end{equation}
the inequality
\begin{equation}
\De\big(|B|^{2p} h_2^q\big)\geq (4q-3p)|B|^{2p+2} h_2^q\label{p2}
\end{equation}
holds. Especially,
\begin{equation}
\De\big(|B|^{p-1} h_2^{\f{p}{2}}\big)\geq \f{3}{2}|B|^{p+1} h_2
^{\f{p}{2}}
\end{equation}
whenever
\begin{equation}
p\geq mn-1.\label{p3}
\end{equation}

Let $\e$ be a smooth function with compact support. Integrating by
parts in conjunction with Young's inequality lead to
\begin{eqnarray}\aligned
\int_M |B|^{2p} h_2^p \eta^{2p}*1
&\leq\f{2}{3}\int_M |B|^{p-1} h_2^{\f{p}{2}}\eta^{2p}\De\big(|B|^{p-1}h_2^{\f{p}{2}}\big)*1\\
&=-\f{2}{3}\int_M \n\big(|B|^{p-1}h_2^{\f{p}{2}}\eta^{2p}\big)\cdot \n\big(|B|^{p-1}h_2^{\f{p}{2}}\big)*1\\
&=-\f{2}{3}\int_M \Big|\n \big(|B|^{p-1}h_2^{\f{p}{2}}\big)\Big|^2
\eta^{2p}*1\\
&\qquad\quad-\f{2}{3}\int_M |B|^{p-1}h_2^{\f{p}{2}}\cdot 2p\eta^{2p-1}\n\eta\cdot \n\big(|B|^{p-1}h_2^{\f{p}{2}}\big)*1\\
&\leq-\f{2}{3}\int_M \Big|\n
\big(|B|^{p-1}h_2^{\f{p}{2}}\big)\Big|^2 \eta^{2p}*1
+\f{2}{3}\int_M \Big|\n \big(|B|^{p-1}h_2^{\f{p}{2}}\big)\Big|^2 \eta^{2p}*1\\
&\hskip1in +\f{2}{3}\int_M p^2 |B|^{2p-2}h_2^p\eta^{2p-2}|\n \eta|^2*1\\
&=\f{2}{3}p^2\int_M |B|^{2p-2}h_2^p\eta^{2p-2}|\n
\eta|^2*1.\endaligned\label{Lp2}
\end{eqnarray}
By H\"{o}lder inequality,
\begin{eqnarray}\aligned
&\int_M |B|^{2p-2} h_2^p \eta^{2p-2}|\n \eta|^2*1\\
&=\int_M |B|^{2p-2} h_2^{p-1}\eta^{2p-2}\cdot h_2|\n \eta|^2*1\\
&\leq\Big(\int_M
|B|^{2p}h_2^p\eta^{2p}*1\Big)^{\f{p-1}{p}}\Big(\int_M h_2^p|\n
\eta|^{2p}*1\Big)^{\f{1}{p}}. \endaligned\label{Lp3}
\end{eqnarray}
By (\ref{Lp2}) and (\ref{Lp3}), we finally arrive at
\begin{equation}
\Big(\int_M |B|^{2p}h_2^p\eta^{2p}*1\Big)^{\f{1}{p}}\leq
\f{2}{3}p^2\Big(\int_M h_2^p|\n
\eta|^{2p}*1\Big)^{\f{1}{p}}.\label{lp4}
\end{equation}

Take $\eta\in C_c^\infty (M_R)$ to be the standard cut-off function
such that $\eta\equiv 1$ in $M_{\th R}$ and  $|\n \eta|\leq C
(1-\th)^{-1}R^{-1}$; then  from (\ref{lp4}) we have the following
estimate.

\begin{thm}\label{p4}
Let $M$ be an $n$-dimensional minimal submanifolds of $\ir{n+m}$. If
the Gauss image of $M_R$ is contained in an open geodesic ball of
radius $\f{\sqrt{2}}{4}\pi$ in $\grs{n}{m}$, then there exists
$C_1=C_1(n, m),$ such that
\begin{equation}\label{Lp4}
\big\||B|^2 h_2\big\|_{L^p(M_{\th R})}\leq
C_2(p)(1-\th)^{-2}R^{-2}\big\|h_2\big\|_{L^p(M_R)}
\end{equation}
whenever $p\geq C_1$ and $\th\in (0, 1)$.
\end{thm}

Furthermore, the mean value inequality for any subharmonic function
on minimal submanifolds in $\ir{m+n}$ (ref. \cite{c-l-y}, \cite{n})
can be applied to yield an estimate of the upper bound of $|B|^2$.

Let $B_R(x)\subset \ir{m+n}$  be a ball of radius $R$ and centered
at $x\in M$. Its restriction on $M$ is denoted by
$$D_R(x)=B_R(x)\cap M$$.

\begin{thm}\label{t2}
Let $x\in M$, $R>0$ such that the image of $D_R(x)$ under the Gauss
map lies in an open geodesic ball of radius $\f{\sqrt{2}}{4}\pi$ in
$\grs{n}{m}$. Then, there exists $C_1=C_1(n,m)$, such that
\begin{equation}\label{ineq14}
|B|^{2p} (x)\leq C(n,p)R^{-(n+2p)}(\sup_{D_R(x)}
h_2)^p\text{Vol}(D_R(x)),
\end{equation}
for arbitrary $p\ge C_1$.
\end{thm}

\begin{proof} Choose $q=p\ge mn-1$ which satisfies (\ref{p1}).
The inequality (\ref{p2}) means that $|B|^{2p}h_2^p$ is a
subharmonic function on the minimal submanifold $M$. By Theorem
\ref{p4} and the mean value inequality,
\begin{eqnarray}\aligned
|B|^{2p} h_2^p(x)&\leq \f{C(n)}{(\f{R}{2})^n}\int_{D_{\f{R}{2}}(x)} |B|^{2p}h_2^p *1\\
&= \f{C(n)}{(\f{R}{2})^n}\big\||B|^2 h_2\big\|^p_{L^p(D_{\f{R}{2}}(x))}\\
&\leq\f{C(n)C_2(p)^p}{(\f{R}{2})^{n+2p}}\|h_2\|_{L^p(D_R(x))}^p\\
&\leq\f{C(n)C_2(p)^p}{(\f{R}{2})^{n+2p}}\big(\sup_{D_R(x)}
h_2\big)^p\ \mbox{Vol}(D_R(x))\endaligned\label{ineq11}
\end{eqnarray}
whenever $p\geq C_1(n,m)$.
\end{proof}

\bigskip\bigskip

\Section{Geometric conclusions}{Geometric conclusions}
\medskip

Let $P_0\in \grs{n}{m}$ be a fixed point which is described by
$$P_0=\ep_1\wedge\cdots\wedge\ep_n,$$
where $\ep_1,\cdots,\ep_n$ are orthonormal vectors in $\ir{m+n}$.
Choose complementary orthonormal vectors
$\ep_{n+1},\cdots,\ep_{n+m}$, such that $\{\ep_1,\cdots,\ep_n,
\ep_{n+1},\cdots,\ep_{n+m}\}$ is an orhtonormal base in $\ir{m+n}$.

Let $p:\ir{m+n}\to \ir{n}$ be the natural projection defined by
$$p(x^1,\cdots,x^n;x^{n+1},\cdots,x^{m+n})=(x^1,\cdots,x^n),$$
which induces a map from $M$ to $\ir{n}.$ It is a smooth map from a
complete manifold to $\ir{n}.$

For any point $x\in M$ choose a local orthonormal tangent frame
field $\{e_1,\cdots, e_n\}$ near $x$. Let $v=v_ie_i\in TM.$ Its
projection
$$p_*v=\left<v_ie_i, \ep_j\right>\ep_j=v_i\left<e_i, \ep_j\right>\ep_j.$$

For any  $P\in\g(M)$,
$$w\mathop{=}\limits^{def.}\left<P, P_0\right>
=\left<e_1\wedge\cdots\wedge e_n, \ep_1\wedge\cdots\wedge
\ep_n\right>=\det W,$$ where $W=(\left<e_i, \ep_j\right>).$   It is
well known that
$$W^TW=O^T\La O,$$
where $O$ is an orthogonal matrix and
$$\La=
\begin{pmatrix} \la_1^2&&0\cr
          &\ddots&\cr
      0&&\la_r^2
  \end{pmatrix},\qquad r=min(m,n),$$
where each $0\le\la_i^2\le 1.$

We now compare the length of any tangent vector $v$ to $M$ with its
projection $p_*v.$
$$|p_*v|^2=\sum_{j=1}^n(v_i\left<e_i, \ep_j\right>)^2=(WV)^TWV,$$
where $V=(v_1, \cdots, v_n)^T$. Hence,
\begin{equation}
|p_*v|^2\ge (\la')^2|v|^2\ge w^2|v|^2\ge w_0^2|v|^2,
\label{di}\end{equation} where $\la'=\min_i\{\la_i\}$ and
$w_0=\inf_Mw$. The induced metric $ds^2$ on $M$ from $\ir{m+n}$ is
complete, so is the homothetic metric $\tilde ds^2=w_0^2ds^2$
whenever $w_0>0$.  (\ref{di}) implies
$$p:(M, \tilde ds^2)\to(\ir{n},\text{canonical metric})$$
increases the distance. It follows that $p$ is a covering map from a
complete manifold into $\ir{n}$ , and a diffeomorphism, since
$\ir{n}$ is simply connected. Hence, the induced Riemannian metric
on $M$ can be expressed as $(\ir{n}, ds^2)$ with
$$ ds^2=g_{ij}dx^i\,dx^j.$$
Furthermore, the immersion $F:M\to \ir{m+n}$ is realized by a graph
$(x, f(x))$ with $f:\ir{n}\to\ir{m}$ and
$$g_{ij}=\de_{ij}+\pd{f^\a}{x^i}\pd{f^\a}{x^j}$$

At each point in $M$ its image $n$-plane $P$ under the Gauss map is
spanned by
$$f_i=\ep_i+\pd{f^\a}{x^i}\ep_\a.$$
It follows that
$$|f_1\wedge\cdots\wedge f_n|^2=\text{det}\left(\de_{ij}
     +\sum_\a\pd{f^\a}{x^i}\pd{f^\a}{x^j}\right)$$
and
$$\sqrt{g}=|f_1\wedge\cdots\wedge f_n|.$$
The $n$-plane $P$ is also spanned by
$$p_i=g^{-\frac 1{2n}}f_i,$$
furthermore,
$$|p_1\wedge\cdots\wedge p_n|=1.$$
We then have
\begin{equation*}\begin{split}
\left<P,P_0\right>&=\text{det}(\left<\ep_i,p_j\right>)\cr
 &= \begin{pmatrix}
g^{-\frac 1{2n}}&&0\\
            &\ddots&\cr
    0&& g^{-\frac 1{2n}}\\
        \end{pmatrix}\\
  &=\frac 1{\sqrt{g}}\ge w_0
\end{split}\end{equation*}
and \begin{equation} \sqrt{g}\le \frac 1{w_0}\label{wup}
\end{equation}

Now,
$$D_R(x)=\big\{(\td x,f(\td x)):\td x\in \Om,
f_1,\cdots,f_m \mbox{ are smooth functions}\ \mbox{on }\Om\big\}$$
where $\Om\subset B_R\subset \ir{n}$. Then (\ref{wup}) implies
\begin{equation}\label{ineq13}
\mbox{Vol}(D_R(x))\le \f{1}{w_0}\cdot \mbox{Vol}(\Om)\leq
\f{1}{w_0}C(n)R^n.
\end{equation}

The previous arguments show that

\begin{pro}\label{vol}
Let $M$ be a complete submanifold in $\ir{m+n}$. If the $w-$function
is bounded below by a positive constant $w_0$. Then $M$ is an entire
graph with Euclidean volume growth. In particular, if the Gauss
image of $M$ is contained in a geodesic ball of radius
$\f{\sqrt{2}}{4}\pi$, then $M$ is an entire graph with Euclidean
volume growth.
\end{pro}
\begin{proof}
Now, we consider the case of the image under the Gauss map $\g$
containing in an open geodesic ball of radius $\frac {\sqrt{2}}{4}
\pi$ and centered at $P_0$. The Jordan angles between $P$ and $P_0$
are
$$\th_i=\cos^{-1}(\la_i),$$
where $\la_i^2$ are eigenvalues of the symmetric matrix $W^TW$
\cite{w}.   We know that
$$w=\prod\cos\th_i.$$
On the other hand, the distance between $P_0$ and $P$ (see
\cite{x2}, pp. 188-194)
$$d(P_0, P)=\sqrt{\sum\th_i^2}$$
which is less than $\frac{\sqrt{2}}{4}\pi$ by the assumption. It
follows that
$$w>w_0=\left(\cos\frac{\sqrt{2}}{4}\pi\right)^r.$$
\end{proof}

Theorem \ref{t1}, Schoen-Simon-Yau's type estimates, and Proposition
\ref{vol} give us the following Bernstein type theorem.
\begin{thm}
Let $M$ be a complete minimal $n$-dimensional submanifold in
$\ir{n+m}$  with $n\leq 6$ and $m\ge 2$. If the Gauss image of $M$
is contained in an open geodesic ball of $\grs{n}{m}$ centered at
$P_0$ and of radius $\f{\sqrt{2}}{4}\pi$, then $M$ has to be an
affine linear subspace.

\end{thm}

\begin{proof} Now we choose
$$p=4+\f{2}{3}+\f{4}{3}\sqrt{1+\f{6}{mn}}>6.$$
Fix $x\in M$ and let $r$ be the Euclidean distance function from $x$
and $M_R=D_R(x)$. Hence, letting $R\to +\infty$ in (\ref{Lp15})
yields
$$\big\||B|\big\|_{L^p(M)}=0.$$
i.e., $|B|^2=0$. $M$ has to be an affine linear subspace.
\end{proof}

Theorem \ref{t2} and Proposition \ref{vol} yield a Bernstein type
result as follows.
\begin{thm}
Let $M$ be a complete minimal $n$-dimensional submanifold in
$\ir{n+m}$. If the Gauss image of $M$ is contained in an open
geodesic ball of $\grs{n}{m}$ centered at $P_0$ and of radius
$\f{\sqrt{2}}{4}\pi$, and $\big(\f{\sqrt{2}}{4}\pi-\rho\circ
\g\big)^{-1}$ has growth
\begin{equation}
\big(\f{\sqrt{2}}{4}\pi-\rho\circ \g\big)^{-1}=o(R),
\end{equation}
where $\rho$ denotes the distance on $\grs{n}{m}$ from $P_0$ and $R$
is the Euclidean distance from any point in $M$. Then $M$ has to be
an affine linear subspace.
\end{thm}

\begin{proof} Now we claim
\begin{equation}
\sec(\sqrt{2}\rho)\leq C\big(\f{\sqrt{2}}{4}\pi-\rho\big)^{-1}
\end{equation}
for a positive constant $C$. It is sufficient to prove the function
$$t\in [0,\f{\sqrt{2}}{4}\pi)\mapsto \sec(\sqrt{2}t)\big(\f{\sqrt{2}}{4}\pi-t\big)$$
is bounded, which follows from
\begin{eqnarray*}
\lim_{t\mapsto(\f{\sqrt{2}}{4}\pi)^-}\sec(\sqrt{2}t)\big(\f{\sqrt{2}}{4}\pi-t\big)=\lim_{t\mapsto(\f{\sqrt{2}}{4}\pi)^-}
\f{\f{\sqrt{2}}{4}\pi-t}{\cos(\sqrt{2}t)}
=\lim_{t\mapsto(\f{\sqrt{2}}{4}\pi)^-}\f{-1}{-\sqrt{2}\sin(\sqrt{2}t)}=\f{\sqrt{2}}{2}.
\end{eqnarray*}
Hence we arrive at the inequality
\begin{equation}
h_2\leq C\left(\f{\sqrt{2}}{4}\pi-\rho\circ\g\right)^{-2}.
\end{equation}
Thus, for any point $x\in M$, by Theorem \ref{t2} and Proposition
\ref{vol}, we have
$$|B|^{2p}(x)\le C(n,
p)R^{-2p}\left(\f{\sqrt{2}}{4}\pi-\rho\circ\g\right)^{-2p}$$ Letting
$R\to +\infty$ in the above inequality forces $|B(x)| =0$.
\end{proof}
From (\ref{Lp13}) it is easy to obtain the following result.

\begin{thm}
Let $M$ be an $n$-dimensional complete minimal submanifolds of
$\ir{n+m}$, if the Gauss image of $M$ is contained in an open
geodesic ball  in $\grs{n}{m}$ of radius $\f{\sqrt{2}}{4}\pi$ and
$M$ has finite total curvature, then $M$ has to be an affine linear
subspace.
\end{thm}

For a minimal $n-$submanifold in $\ir{m+n}$, if its Gauss image is
contained in an open geodesic ball on $\grs{n}{m}$ of radius
$\f{\sqrt{2}}{4}\pi$, there is a positive function
$h_1=\cos(\sqrt{2}\rho\circ\g)$. Then the strong stability
inequality (\ref{lp3}) follows. Besides its key role in S-S-Y's
estimates, there are other applications. We state following results,
whose detail proof can be found in the previous paper of the first
author \cite{x3}.

\begin{thm}\label{van} Let $M$ be a complete minimal $n$-submanifold in $\ir{m+n}.$
If the image under the Gauss map is contained in an open geodesic
ball in $\grs{n}{m}$ of radius $\f{\sqrt{2}}{4}\pi$. Then any
$L^2$-harmonic $1$-form vanishes.
\end{thm}

\begin{thm} Let $M$ be one as in Theorem \ref{van}, $N$ be a manifold
with non-positive sectional curvature. Then any harmonic map $f:M\to
N$ with finite energy has to be constant.
\end{thm}

The argument in \cite{c-s-z} and Theorem \ref{van}  give

\begin{thm}  Let $M$ be one as in Theorem \ref{van}. Then $M$ has only one end.
\end{thm}

\bibliographystyle{amsplain}

\end{document}